\input amstex
\documentstyle{amsppt}

\magnification=1200
\loadmsbm
\loadmsam
\loadeufm
\input amssym
\UseAMSsymbols

\TagsOnRight

\hsize168 true mm
\vsize220 true mm
\voffset=20 true mm
\hoffset=0 true mm
\baselineskip 9 true mm plus0.4 true mm minus0.4 true mm

\def\P{\partial}
\def\l{\ell}

\def\X{\xi}

\topmatter
\title
Towards an Effectivisation of the Riemann Theorem
\endtitle

\author 
S.Natanzon
\endauthor

\affil
Moscow State University, Independent University of Moscow,
Institute for Theoretical and  Experimental Physics
\endaffil

\endtopmatter

\document

\centerline{Abstract}
Let $Q$ be a connected and simply connected domain on the Riemann sphere, not coinciding with the Riemann sphere and with the whole complex plane  $\Bbb C$. Then, according to the Riemann Theorem, there exists a conformal bijection between $Q$ and the exterior of the unit disk. In this paper we find an explicit form of this map for a broad class of domains with analytic boundaries.

2000 Math. Sabj. Class. 30C, 37K.

\subhead
1. Introduction
\endsubhead

Let $Q$ be a connected and simply connected domain on the Riemann sphere, not coinciding with the Riemann sphere and with the whole complex plane  $\Bbb C$. 
Then, according to Riemann Theorem, there exists a conformal bijection between $Q$ and the exterior of the unit disk  $U=\{u\in\bar\Bbb C| |u|>1\}$.  In this paper we find an explicit form of this map for a broad class of domains with analytic boundaries.

Without loss of generally it is possible to assume that $\infty\in Q-\partial Q$ and  $0\in Q_{\circ}$, where $Q_{\circ} = \bar\Bbb C-(Q\cup\partial Q)$ and
$\partial Q$ is the boundary of $Q$. Let $\Bbb H$ be the set of all
such domains. The harmonic moments
$$t_0=\frac 1\pi\int\limits_{Q_{\circ}}dxdy, \quad
t_k=-\frac {1}{\pi k} \int\limits_Q z^{-k}dxdy \quad  (k=1,2...)
\tag 1.1$$
form a coordinate system $t=(t_0,t_1,t_2,...)$ on a considerable part of
$\Bbb H$ [ 2, 14 ]. Note that the evolution of $Q_{\circ}$ as $t_0\to 0$ and $t_i=const$ for $i>1$ describes the evolution of a gas bubble surrounded by luquid [1].

Let $Q^t\in\Bbb H$ be the domain corresponding to
$t$. The conditions $\P_zw_{Q^t}(\infty)\in\Bbb R$,
$\P_zw_{Q^t}(\infty)$ $ > 0$ uniquely define the conformal bijection from $Q^t$ to $U$:

$$w(z,t)=w_{Q^t}(z)=p(t)z+\sum\limits_{j=0}^\infty p_j(t) z^{-j},\quad
\quad  p(t)\in\Bbb R, p(t)>0.\quad \tag 1.2$$

Our goal is to represent $p_j(t)$ as a Taylor series in $t$ and $\bar t$.
 
Put
$$\P_i=\frac{\P}{\P t_i},\quad \bar\P_i=\frac{\P}{\P \bar t_i},
\quad D(z)=\sum_{k\geqslant 1}\frac{z^{-k}}{k}\P_k,\bar D(\bar z)=
\sum_{k\geqslant 1}\frac{\bar z^{-k}}{k}\bar\P_k
\tag 1.3$$

In $\S 1$ we reproduce, following [ 3,15,16 ],  an important result [ 4 ]:
the Riemann maps have a potential  $F : \Bbb H\to \Bbb C$
(termed the tau-function of the curves in [11]) such that

$$w(z,t)=zexp((-\frac 12 \P_0^2-\P_0 D(z))F(t))\tag 1.4$$
and, moreover, $F$ satisfies the differential equations 

$$(z-\X)e^{D(z)D(\X)F}=z e^{-\P_0D(z)F}-\X e^{-\P_0D(\X)F},
\tag 1.5(a)$$
$$(\bar z-\bar\X)e^{\bar D(\bar z)\bar D(\bar\X)F}=
\bar z e^{-\P_0\bar D(\bar z)F}-\bar\X e^{-\P_0\bar D(\bar\X)F},
\tag 1.5(b)$$
$$1- e^{-D(z)\bar D(\bar\X)F}=\frac{1}{z\bar\X} e^{\P_0
(\P_0+D(z)+\bar D(\bar\X))F}. \tag 1.5(c)$$

This system of nonlinear differential equations is well known in
mathematical physics and in theory of integrable systems as
the dispersionless limit of the $2D$ Toda hierarchy [ 5 ].
The solution $F(t)$ satisfies some additional equation,
and it appears in string theory as the
"string solution" [ 6 ]. The string solution of the dispersionless
limit of the $2D$ Toda hierarchy appears also in
matrix models and in some other problems of mathematical
physics [ 7,14 ]. Thus, a description of it has an independent interest.

In $\S 2$ we find, following [ 8 ], recursive formulas for coefficients
of the Taylor series for the potential
$$ F =\sum N(i|i_1..., i_k |\bar i_1..., \bar i_k)
t_0^i t _ {i_1} \dotsb t _ {i_k} \bar t _ {\bar i _ {\bar 1}}
\dotsb\bar t _ {\bar i _ {\bar k}}. \tag 1.6$$
These formulas together with (1.4) allow one to find the Riemann maps via
harmonic moments of the domains.

In $\S 3$ we find, following [ 9 ], a sufficient condition for the convergence
of the Taylor series  (1.6).

I acknowledge useful discussions with S.~P.~Novikov, A.~Marshakov, A.~Zabrodin.
The work was partially supported by the grants RFBR 02-01-22004a and
INTAS-00-0259. The final of this paper was written during the author stay at Max-Planck-Institut f\"{u}r Mathematik in Bonn. I should like
to thank this institution for its support and hospitality.

\subhead
2. The potential for Riemann maps 
\endsubhead

The bijection $w :Q\to U$ generates the Dirichlet Green function on $Q\times Q$,

$$ G_Q(z,\xi)=\log\left\vert\frac{w(z)-w(\xi)}{w(z)
\bar w(\xi)-1}\right\vert. \tag 2.1$$

It solves the Dirichlet problem in $Q$: the formula  
$$ u(z)= -\frac{1}{2\pi}\oint_{\P Q}u_0(\xi)\P_nG_Q(z,\xi)|d\xi| \tag 2.2$$
restores a harmonic function from its boundary values $u_0=u|_{\P D}$.
Here $\P_n$ denotes the derivative along the internal normal
to the boundary $\P Q$ with
respect to the second variable, and $|d\xi|$ is an infinitesimal element of
length along $\P Q$.

The Dirichlet Green function is uniquely determined by the following
properties [ 10 ].

$$G_Q(z,\xi) \text{ is symmetric and harmonic in both
arguments except the line } z=\xi$$
$$ \text{ everywhere in } Q, 
\text{ where } G_Q(z,\xi) =\log|z- \xi|
+O(1) \text{ as } z\to \xi;\tag 2.3(a)$$

$$G_Q(z,\xi)=0\text{ if both variables }z,\xi \text{ belong to the
boundary.} \tag 2.3(b)$$

Denote by $\Bbb H_z$ the set of all domains $Q\in\Bbb H$ containing $z$.
The infinitesimal shift of the boundary $\P Q=\{\xi\}$ on
$-\frac{\varepsilon}{2}\partial_nG_Q(z,\xi)\quad (\varepsilon\to 0)$
generates a vector field $\delta_z$ on $\Bbb H_z$. In particulary,
$$\delta_z(t_k)= -\frac{1}{\varepsilon\pi k}\oint_{\partial Q}\xi^{-k}
(-\frac{\varepsilon}{2}\partial_nG_Q(z,\xi))|d\xi|=\frac{z^k}{k} \quad(k>0)
\quad\quad\quad\quad\delta_z(t_0)=1. \tag 2.4$$

Thus, for any functional $X:\Bbb H_z\to\Bbb C$, we have
$$\delta_z X = \frac{\partial X}{\partial t_0}\delta_z t_0 +
\sum\frac{\partial X}{\partial t_k}\delta_z t_k +
\sum\frac{\partial X}{\partial \bar t_k}\delta_z \bar t_k = \nabla(z)X,
 \tag 2.5$$
where $\nabla(z) =\P_0 + D(z) + \bar D(\bar z)$.

Put $v_0 = \frac {2}{\pi}\int\limits_{Q_{\circ}}\log|z|dxdy$ and
$v_k = \frac {1}{\pi}\int\limits_{Q_{\circ}}z^kdxdy$ for $k>0$.

{\bf Theorem 2.1.} [4,11,16].
{\it The function $F(t)=-\frac 1\pi\iint\limits_{Q_{\circ}^t
Q_{\circ}^t}\log|\eta^{-1}-\nu^{-1}|dx_\eta dy_\eta dx_\nu dy_\nu$ 
(where $Q_{\circ}^t=\bar\Bbb C-Q^t$) satisfies the equations $(1.4)$ and $(1.5)$.
Moreover, $\P_kF = v_k$}.

Proof: It follows from $(2.2)$ and from the definition of $\delta_z$ that
the function
$$\tilde G_Q(z,\xi)=\log\left\vert\frac{1}{z} - \frac{1}{\xi}
\right\vert + \frac 12\delta_z\delta_{\xi}F \tag 2.6$$
satisfies $(2.3)$. Thus, $\tilde G_Q(z,\xi)=
G_Q(z,\xi)$. Using $(2.1)$ and $(2.5)$, we find that

$$\log\left\vert\frac{w(z)-w(\xi)}{w(z)\bar w(\xi)-1}
\right\vert=
\log\left\vert\frac{1}{z} - \frac{1}{\xi}\right\vert + \frac12
\nabla(z)\nabla(\xi)F\tag 2.7$$
This equation implies equations (1.4) and (1.5). Indeed, it implies 

$$h=\log\left\vert\frac{w(z)-w(\xi)}{w(z)\bar w(\xi)-1}
\right\vert^2-
\log\left\vert\frac{1}{z} - \frac{1}{\xi}\right\vert^2-
\nabla(z)\nabla(\xi)F = 0. \tag 2.8$$

Furthermore, we have $h=h_1+h_2$ where
$$h_1=\log\left(\frac{w(z)-w(\xi)}{w(z)\bar w(\xi)-1}\right)-
\log\left(\frac{1}{z} - \frac{1}{\xi}\right) -
\left(\frac12\P_0+D(z)\right)\nabla(\xi)F \tag 2.9$$
is a holomorphic function of $z$ while
$$h_2=\log\left(\frac{\bar w(z)-\bar w(\xi)}{\bar w(z)w(\xi)-1}
\right)-
\log\left(\frac{1}{\bar z} - \frac{1}{\bar \xi}\right) -
\left(\frac12\P_0+\bar D(\bar z)\right)\nabla(\xi)F \tag 2.10$$
is an antiholomorphic function of $z$. Thus, 
$h_1$ is independent of $z$.
Passing on to the limit as $z\to\infty$, we obtain
$$h_1=\log\left(\frac{1}
{\bar w(\xi)}\right)- \log\left(-\frac{1}{\xi}\right) -
\frac{1}{2}\P_0\nabla(\xi)F. \tag 2.11$$
Thus,
$$\log\left(\frac{w(z)-w(\xi)}{z-\xi}
\frac{z\bar w(\xi)}{w(z)\bar w(\xi)-1} \right)=
D(z)\nabla(\xi)F. \tag 2.12$$

Passing on to the limit as $\xi\to\infty$ in $(2.12)$, we obtain
$\log\left(\frac {pz}{w(z)}\right)=D(z)\P_0F$. Substituting $(1.2)$
in $(2.7)$ and passing on to the limit as 
$z\to\infty, \xi\to\infty$,  we obtain 
$\log(p)=-\frac 12\P_0^2F$. The comparison of these formulas yields
$$\log\left(\frac {z}{w(z)}\right)=D(z)\P_0F + \frac 12\P_0^2F, \tag 2.13$$
which is equivalent to $(1.4)$.

Using $(2.13)$, we transform the holomorphic parts of $(2.12)$ into
$$\log\left(\frac{w(z)-w(\xi)}
{z-\xi} \right)=-\frac 12\P_0^2F + D(z)D(\xi)F. \tag 2.14$$
Substituting $(1.4)$ into $(2.14)$ leads to $(1.5(a))$. Changing $(z,\xi)$ to
$(\bar z,\bar\xi)$ in the proof we obtain $(1.5(b))$.

The antiholomorphic part of $(2.12)$ is
$$-\log\left(1 - \frac{1}{w(z)\bar w(\xi)} \right)=
D(z)\bar D(\bar\xi)F. \tag 2.15$$
After substituting $(1.4)$ into $(2.15)$ we obtain $(1.5(c))$.

Moreover, if $z\in\P Q$, then
$$\nabla(z)F=\delta_zF=- \frac{2}{\pi}\int
\limits_{Q^t_{\circ}(t)}\log|\nu^{-1}-z^{-1}|d\nu d{\bar \nu} = v_0 +
2Re\sum_{k\geqslant 1}\frac{v_k}{k}z^{-k}.$$
Thus, $\P_kF= v_k$. $\square$

\subhead
3. The Taylor series for the potential  
\endsubhead

In this section we find, following [ 8 ], recursive formulas for
the coefficients $N$ of the Taylor series
$ F =\sum N(i|i_1..., i_k |\bar i_1..., \bar i_k)
t_0^i t _ {i_1} \dotsb t _ {i_k} \bar t _ {\bar i _ {\bar 1}}
\dotsb\bar t _ {\bar i _ {\bar k}}$
for the potential $F$.

The formulas for $N$ are found according to the following scheme.
At first, using some combina\-torial calculations, we transform
Equation (1.5(a)) into an infinite system of equations
 $$\partial_{i_1}\partial_{i_2}\dotsb\partial_{i_k}F\ =$$ 
$$= \ \sum^\infty _ {m=1} \left (\sum\Sb s_1 +\dotsb +s_m =
i_1 +\dotsb +i_k \\ \ell_1 +\dotsb + \ell_m
=m+k-2 \\ s_j, \ell_j\geqslant 1\endSb
\frac {i_1\dotsb i_k} {s_1\dotsb s_m}
T _ {i_1\dotsb i_k} \pmatrix s_1\dotsb s_m \\
\ell_1\dotsb \ell_m\endpmatrix\partial_0^{\ell_1}\partial_{s_1} 
F\dotsb \partial_0 ^ {\ell_m} \partial _ {s_m} F\right). \tag3.1 $$

During this process we find some recursive formulas for $T$.

Then, using the definition of $F$ as a function on the space
of analytic curves, we find that

$\partial_0 F\left\vert _ {t_0} \right. =-t_0+t_0\ln t_0$ + {\it const} and
$\partial_k F\left\vert _ {t_0} \right. = 0, $ if $k > 0,$

\noindent
where here and later $\left |_{t_0}\right.$ means the restriction
of a function to the straight line
$t_1 =\bar t_1=t_2 =\bar t_2 =\dotsb=0 $.

Because of this formula and from Equation (1.5(c)) it follows that
$$\partial_i\bar\partial_j F\left | _ {t_0} \right. =\cases 0, \
\text {if} \ i\ne j, \\ it^i_0,\ \text {if} \ i=j\ .\endcases $$

Later, using (3.1) and the symmetry of the equations (1.5)
we find, that
$$\partial_i\bar\partial_{i_1}\dotsb\bar\partial_{i_k} 
F\left|_{t_0}\right.=\bar\partial_i\partial_{i_1}\dotsb\partial_{i_k} 
F\left | _ {t=0} \right. =\cases 0, \ \text {if} \ i_1 +\dotsb+i_k\ne i, \\
i_1\dotsb i_k\frac {i!} {(i-k+1)!} t_0^{i-k+1},
\ \text {if}\  i_1 +\dotsb+i_k = i.\endcases$$

This condition and Equation (3.1) give some
recursive formulas for the coefficients $N$. As the final result we get

{\bf Theorem 3.1.} {\it Suppose the formal series $(1.6)$ satisfies in its
 domain of convergence Equations $(1.4)$ and $(1.5)$. Then it has the form,
up to a linear summand,
$$F\ = \ \frac {1} {2}\, t_0^2 \, \log {t_0}\ - \ \frac {3}{4}\, t_0^2
 + \sum\limits_{\Sb k, \bar k, n_r, \bar n_r \geqslant 1 \\
0 < i_1 < \dots < i_k \\ 0 < \bar i_1 < \dots < \bar i _ {\bar k}
\\ i-(n_1+\dots+n_k+\bar n_1+\dots+\bar n_{\bar k})+2\geqslant 0
\endSb}
\frac {i_1 ^ {n_1} \dots i_k^{n_k}} {n_1! \dots n_k!} \ 
\frac {\bar i_1 ^ {\bar n_1} \dots \bar i_{\bar k}^{\bar n _ {\bar k}}}
{\bar n_1! \dots \bar n_{\bar k}!}\cdot$$
$$\cdot N_i^2 \left (\matrix i_1,&\dots,&i_k \\ n_1,&\dots,&n_k \endmatrix 
\left | \matrix \bar i_1,&\dots,&\bar i _ {\bar k} \\ \bar n_1,&\dots,& 
\bar n _ {\bar k} \endmatrix \right.\right) 
t_0^{i-(n_1+\dots+n_k+\bar n_1+\dots+\bar n_{\bar k})
+2} t_{i_1}^{n_1} \dots t_{i_k}^{n_k} \bar 
t_{\bar i_1}^{\bar n_1} \dots \bar t_{\bar i_{\bar k}}
^{\bar n_{\bar k}}$$}

\noindent
{\it where the coefficients $N^2$ can be found by the following recursive
rules:

1) $P_{i, j}(s_1, \dots, s_m) \ =
\ \# \ \{(i_1, \dots, i_m) \ | \ i = i_1 + \dots + i_m, \ 
1 \leqslant i_r \leqslant s_r-1 \},$ where $\# V$ is the cardinality
of the set $V$;

2) $T_{i, j} ^1 (s_1, \dots, s_m) \ = $

$ \sum\limits_{\Sb k \geqslant 1 \\ n_1 + \dots + n_k=m \\ n_r \geqslant 1 \endSb}
\frac {1} {k  n_1! \dots n_k!} \
P_{i, j} \left (\underbrace {s_1+\dots +s_{n_1}}_{n_1},
\dots, \underbrace {s_{n_1+\dots+n_{k-1}+1}+\dots+s_{n_1+\dots+n_k}}_{n_k}
\right);$

\noindent
$T_{i_1, i_2}^2 \left (\matrix s_1,&\dots,&s_m \\ 1,&\dots,&1 \endmatrix
\right) \ = \ T_{i_1, i_2}^1 (s_1, \dots, s_m);$

\noindent
$T_{i_1, \dots, i_k}^2 \left (\matrix s_1,&\dots,&s_m \\ l_1,&\dots,&l_m  
\endmatrix\right) \ = $

\noindent
$ = \ \sum\limits_{\Sb 1 \leqslant i \leqslant j \leqslant m \\ s, l \geqslant 1 \endSb}
l \ T_{s, i_k}^1(s_i, \dots,  s_j)
\ T_{i_1, \dots, i_{k-1}}^2 \left (\matrix s_1,&\dots,&s_{i-1},&s, 
 &s_{j+1},&\dots,&s_m \\ l_1,&\dots,&l_{i-1},&l,&l_{j+1},&\dots,&l_m  
\endmatrix\right)$\ ,

\noindent
where
$s = s_i + \dots + s_j - i_k, \ l = (l_i-1) + \dots + (l_j-1)\ ;$

\noindent\qquad
3) $S_{\bar i_1, \dots, \bar i_{\bar k}} \left (\matrix s_1,&\dots,&s_m  
\\ l_1,&\dots,&l_m \endmatrix\right) \ = $

\vskip 0.4cm
\noindent
$=\ \sum\limits_{\Sb \{\bar i_1^1, \dots, \bar i_1^{n_1} \} \sqcup \dots 
\sqcup \{\bar i_m^1, \dots, \bar i_m^{n_m} \} = 
\{\bar i_1, \dots, \bar i_{\bar k} \} \\ \bar i_r^1 + \dots + 
\bar i_r ^ {n_r} = s_r \\s_r-n_r-\ell_r+1\geqslant 0 \endSb} 
\frac {(s_1-1)!} {(s_1-n_1-l_1+1)! (l_1-1)!} \times \dots \times\
\frac {(s_m-1)!} {(s_m-n_m-l_m+1)! (l_m-1)!}\ ;$

\bigskip\bigskip
\noindent\qquad
4) $N_i^1 (i_1, \dots, i_k | \bar i_1, \dots, \bar i_{\bar k}) \ 
= \ 0, \ \text { if}\ i \ne i_1+\dots+i_k \ \text {or}\ i \ne \bar i_1
+ \dots + \bar i _ {\bar k}\ ;$

\vskip 0.4cm
\noindent
\text {in the other cases}

\vskip 0.4cm
\noindent
$N_i^1 (i | \bar i_1, \dots,  \bar i_{\bar k}) \ = 
\ \frac {(i-1)!} {(i-\bar k+1)!}\ ;\quad$
$N_i^1 (i_1, \dots,  i_k  |  \bar i) \ = \ \frac {(i-1)!} {(i-k+1)!}\ ;$

\vskip 0.4cm
\noindent
$N_i^1 (i_1, \dots, i_k |  \bar i_1, \dots, \bar i _ {\bar k}) \ =$

\vskip 0.4cm
\noindent
$ = \ \sum\limits_{\Sb m \geqslant 1 \\ s_1 + \dots + s_m=i_1 + \dots + i_k \\ 
l_1 + \dots + l_m = m+k-2 \\ s_r, l_r \geqslant 1 \endSb}
(-1) ^ {m+1} \ S _ {\bar i_1, \dots, \bar i _ {\bar k}} \left (\matrix 
s_1,&\dots,&s_m \\ l_1,&\dots,&l_m \endmatrix\right)\ $
$T_{i_1, \dots, i_k} ^2 \left (\matrix s_1,&\dots,&s_m \\
 l_1,&\dots,&l_m \endmatrix\right)\ ,$ if $\ k, \bar k>1\ ;$

\vskip 0.4cm
\noindent
$N_i^2 \left (\matrix i_1,&\dots,&i_k \\ n_1,&\dots,&n_k \endmatrix 
\left | \matrix \bar i_1,&\dots,&\bar i _ {\bar k} \\ 
\bar n_1,&\dots,&\bar n _ {\bar k} \endmatrix \right.\right) \ = $

\vskip 0.2cm
\noindent\qquad
$=\ N_i^1 \left (\underbrace {i_1,  \dots, i_1} _ {n_1}, \dots, 
\underbrace {i_k, \dots, i_k}_ {n_k} \left | \underbrace {\bar i_1, 
\dots, \bar i_1}_{\bar n_1}, \dots, \underbrace {\bar i_{\bar k}, 
\dots, \bar i_{\bar k}}_{\bar n{\bar k}} \right.\right)\ . $
$\square$}

The theorem follows from the lemmas below.

{\bf Lemma 3.1.} {\it The following relations hold
$$z-\sum\limits^\infty_{j=1}\frac {1}{j}z^{-j}\P_1\P_jF=
ze^{-\P_0D(z)}F$$ 
$$\P_1\P_j F=\sum\limits^\infty_{m=1}
\frac {(-1)^{m+1}}{m!}
\sum\limits_{\matrix k_1+\dotsb+k_m=j+1 \\ k_i>0\endmatrix}
\frac{j}{k_1\dotsb k_m}\P_0\P_{k_1} F
\dotsb\P_0\P_{k_m}F.$$}

Proof: According to (1.5(a)), $(z-\X) e^{D(z)D(\X)F}=
(z-\X)(1+(D(z) D(\X)F)+\frac 12(D(z) D(\X)F)^2+\dotsb)=
(z-\X)(1+z^{-1}\X^{-1}\P^2_1 F+z^{-1}\sum\limits_{j=2}^\infty
\frac 1j \X^{-j}\P_1\P_jF+\X^{-1}\sum\limits_{j=2}^\infty
\frac 1j z^{-j}\P_1\P_j F+z^{-2}\X^{-2}f)=(z-\X)+\X^{-1}\P^2_1 F-
z^{-1}\P^2_1 F+\sum\limits_{j=2}^\infty\frac 1j\X^{-j}\P_1\P_j F-
\sum\limits_{j=2}^\infty\frac 1j z^{-j}\P_1\P_j F+
z^{-1}\X^{-1}f.$

On the other hand, according to (1.5(a)) the function
$(z-\X)e^{D(z) D(\X)F}$ is a sum of two functions
$f_1(z)+f_2(\X)$. Thus, 
$f=0$ and $ze^{-\P_0 D(z)F}= z-\sum\limits_{j=1}^\infty\frac 1j
z^{-j}\P_1\P_j F$. Therefore,
$$\sum\limits_{j=1}^\infty
\frac 1j z^{-(j+1)}\P_1\P_j F=1-e^{-\P_0 D(z)F}=
1-\bigl(1+\sum\limits_{m=1}^\infty\frac {(-\P_0 D(z)F)^m}{m!}\bigr)=$$
$$-\sum\limits_{m=1}^\infty\frac{(-1)^m}{m!}
\Bigl(\sum\limits_{k=1}^\infty\frac{z^{-k}}{k}\P_0\P_k F\Bigr)^m=$$
$$=-\sum\limits_{m=1}^\infty\frac{(-1)^m}{m!}
\Bigl(\sum\limits_{n=1}^\infty z^{-n}
\sum\limits_{k_1+\dotsb +k_m=n}\frac{1}{k_1\dotsc k_m}
\P_0\P_{k_1} F\dotsb \P_0 \P_{k_m}F \Bigr)=$$
$$=-\sum\limits_{n=1}^\infty z^{-n}\Bigl(\sum\limits_{m=1}^\infty
\frac{(-1)^m}{m!}\sum\limits_{k_1+\dotsb +k_m=n}
\frac{1}{k_1\dotsc k_m}\P_0\P_{k_1}F\dotsb \P_0\P_{k_m} F\Bigr).$$
Thus,
$$\frac 1j \P_1 \P_j F=-\sum\limits_{m=1}^\infty \frac{(-1)^m}
{m!}\sum\limits_{k_1+\dotsb +k_m=j+1}\frac{1}{k_1\dotsc k_m}
\P_0\P_{k_1}F\dotsb \P_0\P_{k_m} F. \square$$

{\bf Lemma 3.2.} {\it The following relation holds
$$\P_i\P_j F= \sum\limits_{m=1}^\infty\sum\limits_{\matrix s_1+\dotsb +s_m=
j+1 \\ s_i>1\endmatrix}$$
$$\frac{(-1)^{m+1}}{m}\frac{ij}{(s_1-1)\dotsb(s_m-1)}
P_{ij}(s_1-1,...,s_m-1)
\cdot\P_1\P_{s_1-1} F\dotsb\P_1\P_{s_m-1} F,$$
where $P_{ij}(s_1-1,...,s_m-1)$ is the number of representations 
of the form $\{i=i_1+\dotsb+i_m\vert 1\leqslant i_k\leqslant s_k-1,\ k=1,
\dotsb, m\}$ for the number $i$.}

Proof: According to Lemma 3.1 and Equation (1.5(a)),
$(z-\X)e^{D(z) D(\X) F}=$ $$=z-\sum\limits_{j=1}^\infty\frac 1j
z^{-j}\P_1\P_j F-(\X-\sum\limits_{j=1}^\infty\frac 1j
\X^{-j}\P_1\P_j F)=(z-\X)-\sum\limits_{j=1}^\infty\frac 1j
(z^{-j}-\X^{-j})\P_1\P_jF.$$ Thus,
$$e^{D(z) D(\X)}=1+z^{-1}\X^{-1}\sum\limits_{j=1}^\infty\frac 1j
\frac{(z^{-j}-\X^{-j})}{(z^{-1}-\X^{-1})}\P_1\P_j F=$$
$$=1+z^{-1}\X^{-1}\sum\limits_{j=1}^\infty\frac 1j
\big(\sum\limits_{\matrix s+t=j-1 \\ s, t\geqslant 0 \endmatrix}
z^{-s}\X^{-t})\P_1\P_j F=$$
$$=1+\sum\limits_{j=1}^\infty\frac 1j(\sum
\limits_{\matrix s+t=j+1 \\ s,t\geqslant 1\endmatrix}
z^{-s}\X^{-t})\P_1\P_j F.$$ Therefore,
$$D(z)D(\X)F=\sum\limits_{m=1}^\infty\frac{(-1)^{m+1}}{m}
\bigl(\sum\limits_{n=1}^\infty\bigl(\sum\limits_{\matrix s+t=n+1
\\s,t\geqslant 1 \endmatrix}
z^{-s}\X^{-t}\bigr)\frac 1n \P_1\P_n F\bigr)^m=$$
$$=\sum\limits_{j=1}^\infty\frac{(-1)^{m+1}}{m}
\sum\limits_{i,j\geqslant 1} z^{-i}\X^{-j}\cdot$$
$$\Bigl(\sum\limits_{\matrix i_1+\dotsb+i_m=i \\ j_1+\dotsb j_m=j
\\ i_k, j_k\geqslant 1\endmatrix} \frac{1}{i_1+j_1-1}
\P_1\P_{i_1+j_1-1}F\dotsb\frac{1}{i_m+j_m-1}\P_1\P_{i_m+j_m-1}F
\Bigr),$$ that is  
$$\P_i\P_j F= \sum\limits_{m=1}^\infty
\sum\limits_{s_1+\dotsb+s_m=i+j}\frac {(-1)^{m+1}}{m}
\frac{ij}{(s_1-1)\dotsb (s_m-1)}\cdot$$
$$\cdot P_{ij}(s_1-1,...,s_m-1)
\P_1\P_{s_1-1} F\dotsb\P_1\P_{s_m-1} F.\qquad \square$$

{\bf Remark.} The equations
$$\P_i\P_j F= \sum\limits_{m=1}^\infty
\sum\limits_{s_1+\dotsb+s_m=i+j}\frac
{(-1)^{m+1}}{m}\frac{ij}{(s_1-1)\dotsb (s_m-1)}\cdot$$
$$\cdot P_{ij}(s_1-1,...,s_m-1)
\P_1\P_{s_1-1} F\dotsb\P_1\P_{s_m-1} F$$
describe the dispersionless limit of the KP equation.
Another description of this hierarchy is presented in [ 12 ].
A comparison of these descriptions gives some
nontrivial combinatorial identity for $P_{ij}$.

{\bf Lemma 3.3.} {\it The following relation holds
$$\P_i\P_j F= \sum\limits_{m=1}^\infty
\sum\limits_{p_1+\dotsb+p_m=j+i} \frac{ij}{p_1\dotsb p_m}
T_{ij}(p_1\dotsb p_m) \P_0 \P_{p_1} F\dotsb \P_0\P_{p_m} F,$$
where \linebreak $T_{ij}(p_1...p_m)=\sum\limits_
{\matrix n_1+\dotsb+n_k=m \\ n_i>0\endmatrix}
\frac{(-1)^{m+1}}{k}
\frac {1}{n_1!\dotsb n_k!} P_{ij}(p_1+\dotsb+p_{q_1}-1,...,
p_{q_{k-1}+1}+\dotsb+p_{q_k}-1)$, and $q_j=\sum\limits^j_{i=1}n_i$.}

Proof: According to Lemmas 3.1 and 3.2,
$$\P_i\P_j F= \sum\limits_{m=1}^\infty
\sum\limits_{s_1+\dotsb+s_m=j+i} \frac{(-1)^{m+1}}{m}
\frac{ij}{(s_1-1)\dotsb (s_m-1)}
P_{ij}(s_1-1,\dotsb, s_m-1)\cdot$$
$$\cdot\P_1 \P_{s_1-1} F\dotsb \P_1\P_{s_m-1} F)=
\sum\limits_{m=1}^\infty
\sum\limits_{s_1+\dotsb+s_m=j+i}\frac{(-1)^{m+1}}{m}
\frac{ij}{(s_1-1)\dotsb (s_m-1)}\cdot$$ $$\cdot P_{ij}(s_1-1,\dotsb,
s_m-1)\Bigl(\sum^\infty_{n_1=1}\sum\limits_{p_1
+\dotsb+p_{n_1}=s_1}
\frac{(-1)^{n_1+1}}{n_1!}\frac{s_1-1}{p_1\dotsb p_{n_1}}
\P_0 \P_{p_1} F\dotsb \P_0\P_{p_{n_1}}F\Bigl)\dotsb$$
$$\dots\Bigl(\sum^\infty_{n_m=1}
\sum\limits_{p_1+\dotsb +p_{n_m}=s_m}\frac{(-1)^{n_m+1}}{n_m!}
\frac{s_m-1}{p_1\dotsb p_{n_m}}
\P_0\P_{p_1} F\dotsb \P_0\P_{p_{n_m}} F\Bigl)=$$
$$=\sum\limits_{m=1}^\infty
\sum\limits_{p_1+\dotsb+p_m=j+i} \frac{ij}{p_1\dotsb p_m}
T_{ij}(p_1\dotsb p_m) \P_0 \P_{p_1} F\dotsb \P_0\P_{p_m}F. \square$$

By induction, we get from Lemma 3.3

{\bf Lemma 3.4.} {\it The following relation holds
$$\P_{i_1}\P_{i_2}\dotsb \P_{i_k}F=\sum\limits_{m=1}^\infty
\Bigl(\sum\limits_{\matrix s_1+\dotsb+s_m=i_1+\dotsb +i_k \\ \ell_1
+\dotsb +\ell_m=m+k-2 \\ s_j, \ell_j\geqslant 1\endmatrix}
\frac{i_1\dotsb i_k}{s_1\dotsb s_m}\cdot$$
$$\cdot T_{i_1\dotsb i_k}
\pmatrix s_1\dotsb s_m \\ \ell_1\dotsb \ell_m\endpmatrix
\P^{\ell_1}_0\P_{s_1} F\dotsb \P^{\ell_m}_0\P_{s_m} F\Bigr),$$
$$ \text{where}\quad
T_{i_1, i_2}\pmatrix s_1\dotsb s_m \\ \ell_1 \dotsb \ell_m\endpmatrix=
\cases T_{i_1 i_2}(s_1\dotsb s_m),& \ \text{if}\ \ell_1=\dotsb=
\ell_m=1 \\ 0 & \ \text{otherwise}\endcases,$$
$$T_{i_1\dotsb i_k}\pmatrix s_1\dotsb s_m \\ \ell_1\dotsb \ell_m
\endpmatrix=
\sum\limits_{\matrix 1\leqslant i\leqslant j\leqslant m \\ s, \ell>0
\endmatrix} T_{i_1\dotsb i_{k-1}}
\biggl( {\matrix s_1\dotsb s_{i-1} \\ \ell_1\dotsb \ell_{i-1}\endmatrix}
\pmatrix s \\ \ell \endpmatrix {\matrix s_{j+1}\dotsb s_m\\ \ell_{j+1}
\dotsb \ell_m\endmatrix} \biggr)\cdot$$
$$\cdot T_{s,i_k}(s_i, s_{i+1},...,s_j)
\frac{\ell!}{(\ell_i-1)!\dotsb (\ell_j-1)!},$$
$$\text{and} \quad s=s_i+s_{i+1}+\dotsb +s_j-i_k, \quad
\ell=(\ell_i-1)+\dotsb+(\ell_j-1). \square $$}

Define now the Cauchy data for $F$.

{\bf Lemma 3.5.} {\it The following relations hold
$$\P_0 F\vert_{t_0}=-t_0+t_0\ln t_0\quad + \text{const\quad and}
\quad\P_k F\vert_{t_0}=0\quad \text{for}\quad k>0.$$}

Proof: If $t_1=t_2= \dots =0$, then $Q^t_{\circ}=\{w\in \Bbb C||w|<t_0^
{\frac 12}\}$ and thus $w(z,t)\vert_{t_0} = Rz= t_0^{\frac 12}z$.
Now it follows from Theorem 2.1 that
$$\P_0 F\Bigl\vert_{t_0}=\frac 2\pi\int\limits_{\vert z\vert\leqslant R}
\ln\vert z\vert dxdy=\frac 2\pi\int\limits_0^R dr
\int\limits_0^{2\pi} d\varphi\cdot r\cdot\ln\vert r\vert=$$
$$=\frac{2\cdot 2\pi}{2\pi}\int\limits_0^R\ln\vert r\vert
dr^2=2(r^2\ln r\Bigl|_0^R-\int\limits_0^R rdr)=$$
$$=r^2\ln r^2\Bigl\vert_0^R-
2\frac 12 r^2\Bigl\vert_0^R=R^2\ln R^2-R^2=
-t_0+t_0\ln t_0.$$

$\P_kF\vert_{t_0} =
\frac {1}{\pi} \int\limits_{Q_{\circ}^t} z^kdxdy = 0$ for $k>0.$
$\square$

{\bf Lemma 3.6.} {\it The following relations hold
$$\P_i\bar\P_j F\vert_{t_0}=\cases 0 &\ \text{for}\ i\ne j,
\\ i t^i_0 &\ \text{for}\ i=j. \endcases$$}

Proof: It follows from Lemma 3.5 that
$\P_0\P_k F\bigl\vert_{t_0}=0$ for $k>0$ and
$\P^2_0 F\vert_{t_0}=\ln t_0$. Thus,
$$e^{\P_0(\P_0+D(z)+\bar D(\bar\X))F}\Bigl\vert_{t_0}=t_0.$$
Moreover, according to (1.5(c)),
$$1-e^{-D(z)\bar D(\bar\X) F}=z^{-1}\bar\X^{-1} e^{\P_0(\P_0+D(z)+
\bar D(\bar\X))F},$$ and thus
$$-D(z)\bar D(\bar\X) F\bigl\vert_{t_0}=\ln(1-z^{-1}\bar\X^{-1}
t_0)=-\sum\limits_{k=1}^\infty k z^{-k}\X^{-k} t_0^k.$$
Therefore, $\P_i\bar\P_jF\bigl\vert_{t_0}=0$ for $i\ne j$ and 
$\P_i\bar\P_i F\bigl\vert_{t_0}=i t^i_0$ $\square$.

{\bf Lemma 3.7.} {\it The following relations hold
$$\P_{i}\bar\P_{i_1}\dotsb\bar\P_{i_k}F\Bigl\vert_{t_0}=
\bar\P_i\P_{i_1}\dotsb \P_{i_k} F\bigl\vert_{t_0}=
\cases 0, &\ \text{if} \quad i_1
+\dotsb+i_k\ne i \\ i_1\dotsb i_k \frac{i!}{(i-k+1)!} t_0^{i-k+1}, &
\ \text{if}\ i=i_1+\dotsb +i_k \endcases.$$}

Proof: The differentials $\P$ and $\bar\P$ enter Equation
(1.5) in a symmetric way. This gives the first equality.
Moreover, according to Lemmas 3.3 and 3.4, we have
$$\P_{i_1}\P_{i_2}\dotsb \P_{i_k} F=
\frac{i_1\dotsb i_k}{i}T_{i_1\dotsb i_k}{\pmatrix i
\\ k - 1\endpmatrix}\P_0^{k-1}\P_{i} F+$$
$$+\sum\limits_{m=2}^\infty
\Bigl(\sum\limits_{\matrix s_1+\dotsb+s_m=i_1+\dotsb +i_k \\ \ell_1
+\dotsb \ell_m=m+k-2 \\ s_j, t_j\geqslant 1\endmatrix}
\frac{i_1\dotsb i_k}{s_1\dotsb s_m} T_{i_1\dotsb i_k}
\pmatrix s_1\dotsb s_m \\ \ell_1\dotsb \ell_m\endpmatrix\cdot$$
$$\cdot\P^{\ell_1}_0\P_{s_1} F\dotsb \P^{\ell_m}_0\P_{s_m} F
\Bigr)=
\frac{i_1\dotsb i_k}{i}\P_0^{k-1}\P_i F+$$
$$+\sum\limits_{m=2}^\infty
\Bigl(\sum\limits_{\matrix s_1+\dotsb+s_m=i_1+\dotsb +i_k \\ \ell_1
+\dotsb \ell_m=m+k-2 \\ s_j, t_j\geqslant 1\endmatrix}
\frac{i_1\dotsb i_k}{s_1\dotsb s_m} T_{i_1\dotsb i_k}
\pmatrix s_1\dotsb s_m \\ \ell_1\dotsb \ell_m\endpmatrix\cdot$$
$$\cdot\P^{\ell_1}_0\P_{s_1} F\dotsb \P^{\ell_m}_0\P_{s_m} F\Bigr),$$
where $i=i_1+\dotsb+i_k$.
This equality and Lemmas 3.5 and 3.6 together 
give the second equality in the
assertion of Lemma 3.7. $\square$

{\bf Lemma 3.8.} {\it The following relation holds
$$\P_{i_1}\dotsb \P_{i_k}\bar\P_{\bar i_1}\dotsb
\bar\P_{\bar i_{\bar k}} F\Bigl\vert_{t_0}=
\sum_{i=1}^\infty \tilde N_i(i_1\dotsb i_k\bigl\vert
\bar i_1\dotsb \bar i_{\bar k})
t_0^{i-(k+\bar k)+2},\ 
\text{where} \
\tilde N_i\bigl(i_1\dotsb i_k\bigl\vert \bar i_1\dotsb \bar i_{\bar k}
\bigr)=0$$
$ \text{if}\ i-(k+\bar k)+2<0, \ \sum\limits_{j=1}^{k} i_j\ne i$
or $\sum\limits_{j=1}^{\bar k}\bar i_j\ne i$. In opposite case,
$\tilde N_i\bigl(i_1\dotsb i_{i_k}\bigl\vert \bar i_1,\dotsb
\bar i_{\bar k}\bigr)=$ 
$$=\sum\limits_{m=1}^\infty
\Bigl(\sum\limits_{\matrix s_1+\dotsb+s_m=i_1+\dotsb +i_k \\ \ell_1
+\dotsb \ell_m=m+k-2 \\ s_j, \ell_j\geqslant 1\endmatrix}
i_1\dotsb i_k\bar i_k\dotsb \bar i_{\bar k}\cdot$$
$$\cdot T_{i_1\dotsb i_k}
\pmatrix s_1\dotsb s_m \\ \ell_1\dotsb \ell_m\endpmatrix
S_{\bar i_1\dotsb \bar i_{\bar k}}
\pmatrix s_1\dotsb s_m \\ \ell_1\dotsb \ell_m\endpmatrix\Bigr),$$
where  
$$S_{\bar i_1\dotsb \bar i_{\bar k}}
\pmatrix s_1\dotsb s_m \\ \ell_1\dotsb \ell_m\endpmatrix=
\sum\limits_{
\matrix n_1+\dotsb+n _m=k \\ s_i- n_i+1-\ell_i\geqslant 0
\endmatrix} \sum
\frac{(s_1-1)!\dotsb (s_m-1)!}{(s_1- N_1+1-\ell_1)!\dotsb
(s_m-n_m+1-\ell_m)!}$$
and the second summation is carried over all partitions of the set
$\{\bar i_1,...,\bar i_{\bar k}\}$ into subsets
\linebreak $\{j^p_1,...,j^p_{n_p}\}$ $(p=1,...,m)$ such that
$\sum\limits^{n_p}_{\alpha=1}j_\alpha^p=s_p$. $\square$}

Proof: According to Lemma 3.4,
$$\P_{i_1}\dotsb\P_{i_k}
\bar\P_{\bar i_1}\dotsb \bar\P_{\bar i_{\bar k}}F=
\bar\P_{\bar i_{\bar i}}\dotsb \bar\P_{\bar i_{\bar k}}
\Bigl(\sum\limits_{m=1}^\infty
\Bigl(\sum\limits_{\matrix s_1+\dotsb+s_m=i_1+\dotsb +i_k \\ \ell_1
+\dotsb +\ell_m=m+k-2 \\ s_j, \ell_j\geqslant 1\endmatrix}
\frac{i_1\dotsb i_k}{s_1\dotsb s_m}\cdot$$
$$\cdot \ T_{i_1\dotsb i_k}
\pmatrix s_1\dotsb s_m \\ \ell_1\dotsb \ell_m\endpmatrix
\P^{\ell_1}_0\P_{s_1} F\dotsb \P^{\ell_m}_0\P_{s_m} F\Bigr)\Bigl)=$$
$$=\sum\Bigl(\sum\limits_{m=1}^\infty
\Bigl(\sum\limits_{\matrix s_1+\dotsb+s_m=i_1+\dotsb +i_k \\ \ell_1
+\dotsb +\ell_m=m+k-2 \\ s_j, \ell_j\geqslant 1\endmatrix}
\frac{i_1\dotsb i_k}{s_1\dotsb s_m}\ 
T_{i_1\dotsb i_k}\pmatrix s_1\dotsb s_m \\ \ell_1\dotsb
\ell_m\endpmatrix\ \cdot$$
$$\cdot\
\P^{\ell_1}_0 \P_{s_1}
\bar\P_{\bar j_1}^1\dotsb \bar\P_{\bar j_{n_1}}^1
 F\dotsb \P^{\ell_m}_0 \P_{s_m}
\bar\P_{\bar j_1}^m\dotsb \bar\P_{\bar j_{n_m}}^m
 F\Bigr)\Bigl),$$
where the first summation is carried over all partitions of 
$\{\bar i_1,...,\bar i_{\bar k}\}$ into subsets
$\{(j^1_i,...,j^1_{n_1})...(j^m_1,...,$ \linebreak
$j^m_{n_m})\}$.
According to Lemma 3.7, this gives the assertion of
Lemma 3.8 $\square$.

Lemma 3.8 yields

{\bf Lemma 3.9} {\it The following relation holds
$$F=\frac 12 t^2_0\log t_0-\frac 34 t_0^2+
\sum\limits_{\matrix i_1<\dotsb<i_k \\
\bar i_1<\dotsb \bar i_{\bar k}
\endmatrix}\sum_{n_j, \bar n_j=1}^{\infty}
\frac{1}{n_1!\dotsb n_k!\bar n_1!\dotsb\bar n_{\bar k}!}
\tilde N_i\Bigl({\matrix i_1...i_k \\ n_1...n_k\endmatrix}\Bigl\vert
\matrix \bar i_1...\bar i_{\bar k} \\ \bar n_1 ... \bar n_{\bar k}
\endmatrix\Bigr)\cdot$$
$$\cdot t_0^{i-(\sum\limits_{i=1}^k n_i+
\sum\limits_{i=1}^{\bar k}\bar n_i)+2}
t_{i_1}^{n_1}\dotsb t_{i_k}^{n_k} \bar t_{\bar i_1}^{\bar n_1}
\dotsb \bar t_{\bar i_{\bar k}}^{\bar n_{\bar k}},$$
where
$$\tilde N_i\Bigl({\matrix i_1...i_k \\ n_1...n_k\endmatrix}\Bigl\vert
\matrix \bar i_1...\bar i_{\bar k} \\ \bar n_1 ... \bar n_{\bar k}
\endmatrix\Bigr) =0,$$ if
$\ i-(\sum\limits_{i=1}^k n_i+
\sum\limits_{i=1}^{\bar k}\bar n_i)+2<0,$
$i\ne\sum\limits_{j=1}^k n_j i_j$ or $i\ne\sum\limits_{j=1}^{\bar k}
\bar n_j\bar i_j,$ and
$$\tilde N_i\Bigl({\matrix i_1...i_k \\ n_1...n_k\endmatrix}\Bigl\vert
{\matrix \bar i_1...\bar i_{\bar k} \\ \vert n_1 ... \bar n_{\bar k}
\endmatrix}\Bigr) = \tilde N_i(i_1...i_1i_2...i_2...i_k...i_k\bigl\vert
\bar i_1...\bar i_1\bar i_2...\bar i_2...\bar i_{\bar k}...
\bar i_{\bar k}),$$
otherwise; here in the last parentheses each
$i_j$ (respectively, each $\bar i_j$ ) occurs $n_j$ (respectively,
$\bar n_j$) times.}

Theorem 3.1 is equivalent to Lemma 3.9.

Using the algorithm of the present paper,
Yu.~Klimov and A.~Korzh created a computer program for 
calculating the coefficients of the Taylor series $F$.

\subhead
4. Convergence conditions of the Taylor series for the potential  
\endsubhead

The combinatorial coefficients $N^2_i (...)$ have some 
remarkable properties. For example,

\vskip 0.3cm
\proclaim {Theorem 4.1}
$N^2_i\left (\matrix i_1\dotsb i_k \\ n_1\dotsb n_k\endmatrix
\left |\matrix 1 \\ \bar n_1 \endmatrix \right.\right)=
\cases(i-1)!,\ \text {if} \
k=n_1=1, \ i=i_1 =\bar n_1, \\ 0\ \text { otherwise}.
\endcases $
\endproclaim

\demo {Proof} According to our definition,

$$S_{\underbrace{1,  \dots, 1}_{\bar k}} \left(\matrix s_1,&\dots,&s_m
\\ l_1,&\dots,&l_m \endmatrix\right) \ =$$
$$= \ \sum\limits_{\Sb \{\bar i_1^1, \dots, \bar i_1^{n_1}\}
\sqcup \dots \sqcup \{\bar i_m^1, \dots, \bar i_m^{n_m}\} =\\=
\{\underbrace{1, \dots, 1}\} \\ \bar i_r^1 + \dots
+ \bar i_r^{n_r} = s_r \\s_r-n_r-\ell_r+1\geqslant 0 \endSb}
\frac{(s_1-1)!}{(s_1-n_1-l_1+1)!(l_1-1)!} \times \dots \times$$
$$\times\dots\times\frac{(s_m-1)!}{(s_m-n_m-l_m+1)!(l_m-1)!} \ =
\delta_{\ell_1,1}\dots\delta_{\ell_m,1}\frac{\bar k!}{s_1
\dots s_m}.$$

Thus if $k>2$, then

$$N_i^1 (i_1, \dots, i_k  |  1, \dots, 1) \ =$$
$$= \ \sum\limits_{\Sb m \geqslant 1 \\ s_1 + \dots + s_m =
i_1 + \dots + i_k \\ l_1 + \dots + l_m = m+k-2 \\ s_r, l_r \geqslant 1
\endSb}
(-1)^{m+1} \ S_{1, \dots, 1} \left(\matrix s_1,&\dots,&s_m
\\ l_1,&\dots,&l_m \endmatrix\right)\times$$
$$\times \ T_{i_1, \dots, i_k}^2 \left(\matrix s_1,&\dots,&s_m \\
l_1,&\dots,&l_m \endmatrix\right)=0\ .$$
Let now $k=2$ ($i_1, i_2 \geqslant 1$). Then
$$N_i^1 (i_1, i_2  | \underbrace{ 1, \dots, 1}_{\bar k}) \ =$$
$$= \ \sum\limits_{\Sb m \geqslant 1 \\ s_1 + \dots + s_m = i_1 + i_2
\\ s_r \geqslant 1 \endSb}
(-1)^{m+1} \ S_{\underbrace{ 1, \dots, 1}_{\bar k}}
\left(\matrix s_1,&\dots,&s_m \\ 1,&\dots,&1 \endmatrix\right)\times$$
$$\times\ T_{i_1, i_2}^2 \left(\matrix s_1,&\dots,&s_m
\\ 1,&\dots,&1 \endmatrix\right) \ =$$
$$= \ \sum\limits_{\Sb m \geqslant 1 \\ s_1 + \dots + s_m =
i_1 + i_2 \\ s_r \geqslant 1 \endSb}
(-1)^{m+1}  \frac{\bar k!}{s_1  \dots  s_m}
\ T_{i_1, i_2}^2 \left(\matrix s_1,&\dots,&s_m \\ 1,&\dots,&1
\endmatrix\right) \ =$$
$$= \ \sum\limits_{\Sb m \geqslant 1 \\ s_1 + \dots + s_m = i_1 + i_2
\\ s_r \geqslant 1 \endSb}
(-1)^{m+1}   \frac{\bar k!}{s_1  \dots  s_m}
\ T_{i_1, i_2}^1 (s_1,  \dots,  s_m) \ =$$
$$= \ \sum\limits_{\Sb m \geqslant 1 \\ s_1 + \dots + s_m = i_1 + i_2
\\ s_r \geqslant 1 \endSb}
(-1)^{m+1}  \frac{\bar k!}{s_1  \dots  s_m}  
 \sum\limits_{\Sb k \geqslant 1 \\ n_1 + \dots + n_k=m
\\ n_r \geqslant 1 \endSb}
\frac {1}{k  n_1!  \dots  n_k!} \times$$
$$\times \ P_{i, j} \left(\underbrace{s_1 + \dots + s_{n_1}}_{n_1},
\dots, \underbrace{s_{n_1 + \dots + n_{k-1} + 1} + \dots +
s_{n_1 + \dots + n_k}}_{n_k}\right) \ =$$

$$= \ \sum\limits_{\Sb m \geqslant 1 \\ s_1 + \dots + s_m = i_1 + i_2
\\ k \geqslant 1 \\ n_1 + \dots + n_k=m \endSb}
\frac {(-1)^{m+1}}{k  n_1!  \dots  n_k!}  \frac{\bar k!}{s_1
 \dots  s_m}\times$$
$$\times\ P_{i, j} \left(\underbrace{s_1 + \dots + s_{n_1}}_{n_1},
\dots, \underbrace{s_{n_1 + \dots + n_{k-1} + 1} + \dots +
s_{n_1 + \dots + n_k}}_{n_k}\right) \ =$$
$$= \ \sum\limits_{\Sb k \geqslant 1 \\ \tilde s_1 + \dots + \tilde s_k
= i_1 + i_2 \endSb}
\ P_{i, j} (\tilde s_1, \dots,  \tilde s_k)\times$$
$$\times\ \sum\limits_{\Sb n_r \geqslant 1,\ n_1+\dots+n_k = m
\\ \underbrace{s_{n_1 + \dots + n_{r-1} + 1} + \dots +
s_{n_1 + \dots + n_r}}_{n_r} = \tilde s_r \endSb}
\frac {(-1)^{m+1}}{k  n_1!  \dots  n_k!} 
\frac{\bar k!}{s_1  \dots  s_m} \ =$$
$$= \ \sum\limits_{\Sb k \geqslant 1 \\ \tilde s_1 + \dots +
\tilde s_k = i_1 + i_2 \endSb}
-\frac{\bar k!}{k}\ P_{i, j} (\tilde s_1, \dots,  \tilde s_k)\times$$
$$\times\prod\limits_{1 \leqslant r \leqslant k} \sum\limits_{\Sb n_r \geqslant 1
\\ s_1 + \dots + s_{n_r} = \tilde s_r \endSb}
\frac {(-1)^{n_r}}{n_r!  s_1  \dots  s_{n_r}}\ .$$

In addition, if $s>1$, then 

$$\sum\limits_{\Sb n \geqslant 1 \\ s_1 + \dots + s_{n} =
s \endSb} \frac {(-1)^{n}}{n!  s_1  \dots  s_{n}} \ =
\ \frac{1}{s!}\frac{\partial^s}{\partial x^s}\
\sum\limits_{n \geqslant 1} \frac {(-1)^{n}}{n!} \left.
\left(x+\frac{x^2}{2}+\frac{x^3}{3}+\dots \right)^n \right|_{x=0}\ =$$
$$=\ \frac{1}{s!}\frac{\partial^s}{\partial x^s}\
\sum\limits_{n \geqslant 1} \frac {(-1)^{n}}{n!} \left. (-\log(1-x))^n
\right|_{x=0} \ = $$
$$= \ \frac{1}{s!}\frac{\partial^s}{\partial x^s}\
\left. \sum\limits_{n \geqslant 1} \frac {(\log(1-x))^n}{n!}
\right|_{x-0} \ =
\ \frac{1}{s!}\frac{\partial^s}{\partial x^s}\
\left. ( \exp(\log(1-x))-1) \right|_{x=0} \ = \ 0.$$
Thus 
$$N_i^1 (i_1, i_2  | \underbrace{ 1, \dots, 1}_{\bar k}) \ =
\ \sum\limits_{\Sb k \geqslant 1 \\ \tilde s_1 + \dots + \tilde s_k
= i_1 + i_2 \endSb}
-\frac{\bar k!}{k}\ P_{i, j} (\tilde s_1, \dots,  \tilde s_k)\times$$
$$\times\ \prod\limits_{1 \leqslant r \leqslant k} \sum\limits_{\Sb n_r \geqslant 1
\\ s_1 + \dots + s_{n_r} = \tilde s_r \endSb}
\frac {(-1)^{n_r}}{n_r!  s_1  \dots  s_{n_r}}=0\ ,$$
since $P_{i, j} (1, \dots,  1) \ = 0.$

If $k=1$, then it follows from our definition that

$$N^2_i\left (\matrix i_1 \\ n_1\endmatrix
\left |\matrix 1 \\ \bar n_1 \endmatrix \right.\right)
=\cases(i-1)!,\ \text {if} \
n_1=1, \ i=i_1 =\bar n_1, \\ 0\ \text { otherwise}.
\endcases\ .\ \qed$$
\enddemo

If $t_i, \bar t_i=0$ for $i>2$ , then $\P Q$ is an ellipse. In this case
Theorem 3.1 gives
$$F=-\frac 34t^2_0+\frac 12 t^2_0\log\Big(\frac{t_0}{1-4|t_2|^2}\Big)
+\frac{t_0}{1-4|t_2|^2}\big(|t_1|^2+t^2_1\bar t_2+\bar t^2_1t_2\big).$$
This formula was obtained first in [ 4 ] by using formulas
for conformal maps from an ellipse to the circle.

The recursive formulas for coefficients of the Taylor series $F$ give a
possibility to estimate the coefficients and to find sufficient convergence
conditions for $F$ provided $t_i, \bar t_i=0 $ for $i>n$.

\proclaim
{Theorem 4.2} Let $\tilde t = (t_0,t_1,\bar t_1,t_2,\bar t_2,...)$ be such that
$t_i, \bar t_i=0 $ for $i>n, 0<t_0<1$ and $|t_i|,\ | \bar t_i |
\leqslant (4n^3 2^n e^n)^{-1}$. Then the series $F(\tilde t)$ converges.
\endproclaim

The proof is based on a sequence of estimations of all values
used in the definition of $N^2$. Here are the estimates:

\proclaim
{Lemma 4.1}  Let $i + j \ = \ s_1 + \dots + s_m$. Then
$P _ {ij} (s_1..., s_m) \leqslant\min (C^{m-1}_{i-1}, C^{m-1}_{j-1})$.
\endproclaim
\demo {Proof}
$$P_{j, i}(s_1, \dots, s_m) \ =P_{i, j}(s_1, \dots, s_m) \ 
= \ \# \{ (i_1, \dots, i_m) \ | \ i \ =$$
$$=\ i_1 + \dots + i_m,\ 1 \leqslant i_r \leqslant s_r-1 \} \ \leqslant$$
$$\leqslant \ \# \{ (i_1, \dots, i_m) | i =
i_1 + \dots + i_m, 1 \leqslant i_r \} \ = \ C_{i-1}^{m-1}\ . \ \qed $$
\enddemo

\proclaim
{Lemma 4.2}  Let $i + j \ = \ s_1 + \dots + s_m$. Then
$T^1_{ij}(s_1,...,s_m)\leqslant\frac{\ell^{m-1}}{m!},\
\text {where} \ \ell =\min (i, j) $.
\endproclaim

\demo{Proof}
$$T_{j, i}^1 (s_1, \dots, s_m) \ =T_{i, j}^1 (s_1, \dots, s_m) \
=\ \sum\limits_{\Sb k \geqslant 1 \\ n_1 + \dots + n_k = m \\ n_r
\geqslant 1 \endSb}\frac {1}{k  n_1!  \dots  n_k!}\ \times$$
$$\times \ P_{i, j}\left(\underbrace{s_1 + \dots + s_{n_1}}_{n_1},
\dots, \underbrace{s_{n_1 + \dots + n_{k-1} + 1} + \dots + s_{n_1 +
\dots + n_k}}_{n_k}\right) \ \leqslant$$
$$\leqslant \ \sum\limits_{\Sb k \geqslant 1 \\ n_1
+ \dots + n_k = m \\ n_r \geqslant 1
\endSb} \frac {C_{i-1}^{k-1}}{k  n_1!  \dots  n_k!} \ =$$
$$= \ \frac{1}{m!}\ \frac{\partial^m}{\partial x^m}\
\frac{1}{i}  \sum\limits_{k \geqslant 1} C_{i}^k  \sum
\limits_{n_1 + \dots + n_k = m} \frac {1}{n_1! \dots n_k!} \ =\ $$
$$\ \frac{1}{m!}\ \frac{\partial^m}{\partial x^m}\ \frac{1}{i}
\sum\limits_{k \geqslant 1}  C_{i}^k  \left.
\left(x + \frac{x^2}{2!} + \dots \right)^k \right|_{x=0} \ =$$
$$=\ \frac{1}{m!}\ \frac{\partial^m}{\partial x^m}\
\frac{1}{i}  \left. \left( \sum\limits_{k \geqslant 1}
C_{i}^k (e^x-1)^k \right) \right|_{x=0} \ =$$
$$= \ \frac{1}{m!}\ \frac{\partial^m}{\partial x^m}\
\frac{1}{i}  \left. ((1+(e^x-1))^i-1) \right|_{x=0} \
= \ \frac{1}{m!}\frac{\partial^m}{\partial x^m}\
\frac{1}{i}  \left. (e^{ix}-1) \right|_{x=0} \ =
\ \frac{i^{m-1}}{m!}\ . \ \qed$$
\enddemo

\proclaim
{Lemma 4.3} Let $\ i_1 + \dots + i_k\ = \ s_1 + \dots + s_m$
and $ (\ell_1-1) + \dotsb + (\ell_m-1) =k-2.$ Then
$T^2_{i_1..., i_k} \pmatrix s_1 \dotsb s_m \\ \ell_1
\dotsb\ell_m\endpmatrix\ \leqslant\frac{I^{m-1}(k-1)^m(k-2)!}{m!}\ ,$ 
where $I =\max (i_r) $.
\endproclaim

\demo {Proof}
We use induction by $k$. If $k=2$, then
$$T_{i_1, i_2}^2 \left(\matrix s_1,&\dots,&s_m \\ 1,&\dots,&1
\endmatrix\right) \ = \
T_{i_1, i_2}^1 (s_1, \dots, s_m) \ \leqslant$$
$$\leqslant \ \frac{I^{m-1}}{m!} \
= \ \left.\frac{I^{m-1}(k-1)^{m}(k-2)!}{m!} \right|_{k \ =\ 2}\ .$$

Let $k>2$. Note firstly that if $l = (l_i-1) + \dots + (l_j-1)$, then 
\nopagebreak
$$\sum\limits_{1 \leqslant i \leqslant j \leqslant m,\ j-i=d}
\frac{l}{(k-2)(d+1)} \ = \ \frac{\sum\limits_{1 \leqslant i
\leqslant j \leqslant m,\ j-i=d}  (l_i-1)  +
\dots  +  (l_j-1)}{(k-2)(d+1)} \ \leqslant$$
$$\leqslant \ \frac{\sum\limits_{1 \leqslant t \leqslant m}
(l_t-1) (d+1)}{(k-2)(d+1)} \ = \ 1.$$
Thus
\nopagebreak
$$T_{i_1, \dots, i_k}^2 \left(\matrix s_1,&\dots,&s_m \\
l_1,&\dots,&l_m \endmatrix\right) \ =
\ \sum\limits_{\Sb 1 \leqslant i \leqslant j \leqslant m \\ s,\ l
\geqslant 1 \endSb} l \ T_{s, i_k}^1 (s_i, \dots, s_j) \ \times$$
$$\times \ T_{i_1, \dots, i_{k-1}}^2 \left(\matrix
s_1,&\dots,&s_{i-1},&s,&s_{j+1},&\dots,&s_m \\ l_1,&
\dots,&l_{i-1},&l,&l_{j+1},&\dots,&l_m \endmatrix\right) \ \leqslant$$

$$\leqslant \ \sum\limits_{1 \leqslant i \leqslant j \leqslant m} l
\frac{I^{j-i}}{(j-i+1)!}  \frac{I^{m-(j-i+1)+1-1} 
(k-2)^{m-(j-i+1)+1}  (k-3)!}{(m-(j-i+1)+1)!} \ =$$
$$= \ I^{m-1}  (k-3)!  \sum\limits_{1 \leqslant i \leqslant j \leqslant m} l
 \frac{(k-2)^{m-(j-i)}}{(j-i+1)!  (m-(j-i))!} \ = $$
$$= \ \frac{I^{m-1}  (k-2)!}{m!}  \sum\limits_{1 \leqslant i
\leqslant j \leqslant m,\ d=j-i}  \frac{l}{(k-2)(d+1)}  \frac{m!}{d!
 (m-d)!}  (k-2)^{m-d} \ =$$
$$= \ \frac{I^{m-1}  (k-2)!}{m!}  \sum\limits_{d=0}^{m-1}
\frac{m!}{d!  (m-d)!}  (k-2)^{m-d}  \left(
\sum\limits_{1 \leqslant i \leqslant j \leqslant m,\ j-i=d} 
\frac{l}{(k-2)(d+1)} \right) \ \leqslant$$
\nopagebreak
$$\leqslant \ \frac{I^{m-1}  (k-2)!}{m!}  \sum\limits_{d=0}^{m-1}
\frac{m!}{d!  (m-d)!}  (k-2)^{m-d} \ \leqslant
\ \frac{I^{m-1}  (k-2)!}{m!}  \sum\limits_{t=0, t=m-d}^{m}
C_m^{t}  (k-2)^{t} \ =$$
$$= \ \frac{I^{m-1}  (k-2)!}{m!}  ((k-2)+1)^m \ =
\ \frac{I^{m-1}  (k-1)^m  (k-2)!}{m!}\ .\qed $$
\enddemo

\proclaim
{Lemma 4.4} Let $\bar I = max (\bar i_r)$ and 
$$\tilde S_{\bar i_1, \dots, \bar i_{\bar k}} (m, k) \ =
\sum\limits_{\Sb \{\bar i_1^1, \dots, \bar i_1^{n_1}\}
\sqcup \dots \sqcup \{\bar i_m^1, \dots, \bar i_m^{n_m}\} =\\
= \{\bar i_1, \dots, \bar i_{\bar k}\} \\ (l_1-1) + \dots + (l_m-1)
= k-2 \\ n_r,\ l_r \geqslant 1 \\ s_r=\bar i_r^1 + \dots + \bar i_r^{n_r}
\\ s_r-n_r-\ell_r+1\geqslant 0 \endSb}
\frac{(s_1-1)!}{(s_1-n_1-l_1+1)!  (l_1-1)!} \times \dots \times$$
$$\times \ \frac{(s_m-1)!}{(s_m-n_m-l_m+1)!  (l_m-1)!}\ .$$

Then $\tilde S_{\bar i_1, \dots, \bar i_{\bar k}}
(m, k) \ \leqslant \ m  (\bar k-1)!  C_{\bar I \bar k - \bar k}^{k-2}
 C_{\bar I \bar k}^{\bar k-m}$.
\endproclaim
 
\demo{Proof} We use the equality
$$\sum\limits_{\tilde n_1 + \dots + \tilde n_m = \bar k-m}
C_{\bar I \tilde n_1 + \bar I}^{\tilde n_1} \times \dots \times
C_{\bar I \tilde n_m + \bar I}^{\tilde n_m}  \frac{\bar I}
{\bar I \tilde n_1 + \bar I} \times \dots \times \frac{\bar I}
{\bar I \tilde n_m + \bar I} \ =$$
$$=\ C_{\bar I (\bar k-m) + m \bar I}^{\bar k-m} 
\frac{m\bar I}{\bar I (\bar k-m) + m\bar I}$$
which follows from [ 13 ]. Then
$$\tilde S_{\bar i_1, \dots, \bar i_{\bar k}} (m, k) \ =$$
$$= \ \sum\limits_{\Sb \{\bar i_1^1, \dots, \bar i_1^{n_1}\}
\sqcup \dots \sqcup \{\bar i_m^1, \dots, \bar i_m^{n_m}\}
=\\= \{\bar i_1, \dots, \bar i_{\bar k}\} \\
(l_1-1) + \dots + (l_m-1) = k-2 \\ n_r,\ l_r \geqslant 1 \\ s_r
= \bar i_r^1 + \dots + \bar i_r^{n_r}
\\ s_r-n_r-\ell_r+1\geqslant 0 \endSb}
\frac{(s_1-1)!}{(s_1-n_1-l_1+1)!  (l_1-1)!} \times \dots \times $$
$$\times \frac{(s_m-1)!}{(s_m-n_m-l_m+1)!  (l_m-1)!} \ =$$
$$= \ \sum\limits_{\Sb \{\bar i_1^1, \dots, \bar i_1^{n_1}\}
\sqcup \dots \sqcup \{\bar i_m^1, \dots, \bar i_m^{n_m}\}
=\\= \{\bar i_1, \dots, \bar i_{\bar k}\} \\
s_r = \bar i_r^1 + \dots + \bar i_r^{n_r}
\geqslant n_r \geqslant 1 \endSb}
\frac{(s_1-1)!}{(s_1-n_1)!} \times \dots \times
\frac{(s_m-1)!}{(s_m-n_m)!}  \times$$
$$\times  \sum\limits_{\Sb (l_1-1) + \dots + (l_m-1) = k-2
\\ s_r-n_r-\ell_r+1\geqslant 0 \endSb}
\frac{(s_1-n_1)!}{(s_1-n_1-l_1+1)!  (l_1-1)!} \times \dots\times$$
$$\times\ \frac{(s_m-n_m)!}{(s_m-n_m-l_m+1)!  (l_m-1)!} \ =$$
$$= \ \sum\limits_{\Sb \{\bar i_1^1, \dots, \bar i_1^{n_1}\}
\sqcup \dots \sqcup \{\bar i_m^1, \dots, \bar i_m^{n_m}\}
=\\= \{\bar i_1, \dots, \bar i_{\bar k}\}
\\ s_r = \bar i_r^1 + \dots + \bar i_r^{n_r}
\geqslant n_r \geqslant 1 \endSb}
\frac{(s_1-1)!}{(s_1-n_1)!} \times \dots \times$$
$$\times \frac{(s_m-1)!}{(s_m-n_m)!}
\sum\limits_{\tilde l_1 + \dots + \tilde l_m = k-2}
C_{s_1-n_1}^{\tilde l_1}  \dots  C_{s_m-n_m}^{\tilde l_m} \ =$$
$$= \ \sum\limits_{\Sb \{\bar i_1^1, \dots, \bar i_1^{n_1}\}
\sqcup \dots \sqcup \{\bar i_m^1, \dots, \bar i_m^{n_m}\}
=\\= \{\bar i_1, \dots, \bar i_{\bar k}\}
\\ s_r = \bar i_r^1 + \dots + \bar i_r^{n_r}
\geqslant n_r \geqslant 1  \endSb}
\frac{(s_1-1)!}{(s_1-n_1)!} \times \dots \times$$
$$\times\ \frac{(s_m-1)!}{(s_m-n_m)!}  C_{(s_1+\dots+s_m)-
(n_1+\dots+n_m)}^{k-2} \ =$$
$$= \ C_{\bar i_1 + \dots + \bar i_{\bar k} - \bar k}^{k-2} 
\sum\limits_{\Sb \{\bar i_1^1, \dots, \bar i_1^{n_1}\}
\sqcup \dots \sqcup \{\bar i_m^1, \dots, \bar i_m^{n_m}\}
=\\= \{\bar i_1, \dots, \bar i_{\bar k}\}
\\ s_r = \bar i_r^1 + \dots + \bar i_r^{n_r}
\geqslant n_r \geqslant 1 \endSb}
\frac{(s_1-1)!}{(s_1-n_1)!} \times \dots \times
\frac{(s_m-1)!}{(s_m-n_m)!}\ \leqslant\ $$
$$\leqslant\ \tilde S_{\underbrace{\bar I, \dots, \bar I}_{\bar k}} (m, k) \
= \ C_{\bar I \bar k - \bar k}^{k-2}  \sum\limits_{\Sb
\{\bar i_1^1, \dots, \bar i_1^{n_1}\} \sqcup \dots \sqcup
\\ \sqcup \{\bar i_m^1, \dots, \bar i_m^{n_m}\}
=\\= \{\bar I, \dots, \bar I\},\ n_r \geqslant 1 \\ s_r
= \bar i_r^1 + \dots + \bar i_r^{n_r} = I n_r\endSb}
\frac{(s_1-1)!}{(s_1-n_1)!} \times \dots \times
\frac{(s_m-1)!}{(s_m-n_m)!} \ =$$
$$= \ C_{\bar I \bar k - \bar k}^{k-2}  \sum\limits_{\Sb
n_1 + \dots + n_m = \bar k \\ n_r \geqslant 1 \endSb}
\frac{\bar k!}{n_1!  \dots  n_m!} 
\frac{(\bar I n_1 - 1)!}{(\bar I n_1 - n_1)!} \times \dots \times
\frac{(\bar I n_m - 1)!}{(\bar I n_m - n_m)!} \ =$$
$$= \ \bar k!  C_{\bar I \bar k - \bar k}^{k-2} 
\sum\limits_{\Sb \tilde n_1 + \dots + \tilde n_m = \bar k - m
\\ \tilde n_r = n_r -1 \geqslant 0 \endSb}
\frac{(\bar I \tilde n_1 + \bar I - 1)!}{(\tilde n_1+1)!
(\bar I \tilde n_1 + \bar I - \tilde n_1 -1)!} \times \dots\times$$
$$\times\ \frac{(\bar I \tilde n_m + \bar I - 1)!}{(\tilde n_m+1)!
 (\bar I \tilde n_m + \bar I - \tilde n_m -1)!} \ =$$
$$= \ \bar k!  C_{\bar I \bar k - \bar k}^{k-2} 
\sum\limits_{\tilde n_1 + \dots + \tilde n_m =
\bar k - m} \frac{\bar I \tilde n_1 + \bar I -
\tilde n_1}{(\bar I \tilde n_1 + \bar I)  (\tilde n_1+1)}
 \frac{(\bar I \tilde n_1 + \bar I)!}{\tilde n_1! 
(\bar I \tilde n_1 + \bar I - \tilde n_1)!} \times \dots \times$$
$$\times\ \frac{\bar I \tilde n_m + \bar I - \tilde n_m}{(\bar I
\tilde n_m + \bar I)  (\tilde n_m+1)}  \frac{(\bar I
\tilde n_m + \bar I)!}{\tilde n_m!  (\bar I \tilde n_m +
\bar I - \tilde n_m)!} \ \leqslant$$
$$\leqslant \ \bar k!  C_{\bar I \bar k - \bar k}^{k-2} 
\sum\limits_{\tilde n_1 + \dots + \tilde n_m = \bar k - m}
\frac{\bar I}{\bar I \tilde n_1 + \bar I}  C_{\bar I \tilde n_1
+ \bar I}^{\tilde n_1} \times \dots \times
\frac{\bar I}{\bar I \tilde n_m + \bar I}  C_{\bar I \tilde n_m
+ \bar I}^{\tilde n_m} \ =$$
$$= \ \bar k!  C_{\bar I \bar k - \bar k}^{k-2} 
C_{\bar I (\bar k-m) + m \bar I}^{\bar k-m} 
\frac{m\bar I}{\bar I (\bar k-m) + m\bar I} \ \leqslant
\ m  (\bar k-1)!  C_{\bar I \bar k - \bar k}^{k-2} 
C_{\bar I \bar k}^{\bar k-m}\ . \ \qed$$
\enddemo

\proclaim
{Lemma 4.5}  $N^1_i (i_1..., i_k |\bar i_1..., \bar i _ {\bar k}) \ \leqslant
\ (k-1)! (\bar k-1)! e^{I(k-1)} 2^{\bar I \bar k - \bar k}  2^{\bar I
\bar k}.$
\endproclaim

\demo {Proof}
$$N_i^1 (i_1, \dots, i_k  |  \bar i_1, \dots, \bar i_{\bar k}) \ = $$
$$= \ \sum\limits_{\Sb m \geqslant 1 \\ s_1 + \dots + s_m = i_1 + \dots +
i_k \\ l_1 + \dots + l_m = m+k-2 \\ s_r, l_r \geqslant 1 \endSb}
(-1)^{m+1} \ S_{\bar i_1, \dots, \bar i_{\bar k}} \left(
\matrix s_1,&\dots,&s_m \\ l_1,&\dots,&l_m \endmatrix\right)  \times$$
$$\times \ T_{i_1, \dots, i_k}^2 \left(\matrix
s_1,&\dots,&s_m \\ l_1,&\dots,&l_m \endmatrix\right) \ \leqslant$$
$$\leqslant\ \sum\limits_{\Sb m \geqslant 1\\s_1
+\dots+s_m=i_1+\dots+i_k\\l_1+\dots+l_m
=m+k-2 \\ s_r,l_r \geqslant 1 \endSb}
\frac{I^{m-1} (k-1)^{m}  (k-2)!}{m!}  \times$$
$$\times  \sum\limits_{\Sb \{\bar i_1^1, \dots,
\bar i_1^{n_1}\} \sqcup \dots \sqcup \{\bar i_m^1, \dots,
\bar i_m^{n_m}\} =\\= \{\bar i_1, \dots, \bar i_{\bar k}\} \\
\bar i_r^1 + \dots + \bar i_r^{n_r} = s_r
\\s_r-n_r-\ell_r+1\geqslant 0 \endSb}
\frac{(s_1-1)!}{(s_1-n_1-l_1+1)!  (l_1-1)!} \times \dots \times$$
$$\times \ \frac{(s_m-1)!}{(s_m-n_m-l_m+1)!  (l_m-1)!} \ =$$
$$= \ \sum\limits_{m \geqslant 1}
\frac{I^{m-1} (k-1)^{m-1} (k-2)!}{m!}  \times$$
$$\times  \sum\limits_{\Sb \{\bar i_1^1, \dots,
\bar i_1^{n_1}\} \sqcup \dots \sqcup \{\bar i_m^1,
\dots, \bar i_m^{n_m}\} =\\= \{\bar i_1, \dots,
\bar i_{\bar k}\} \\ (l_1-1) + \dots + (l_m-1) =
k-2 \\ n_r,\ l_r \geqslant 1 \\ s_r = \bar i_r^1 + \dots
+ \bar i_r^{n_r} \\s_r-n_r-\ell_r+1\geqslant 0 \endSb}
\frac{(s_1-1)!}{(s_1-n_1-l_1+1)! (l_1-1)!} \times \dots\times$$
$$\times\frac{(s_m-1)!}{(s_m-n_m-l_m+1)! (l_m-1)!} \ =$$
$$= \ \sum\limits_{m \geqslant 1}
\frac{I^{m-1} (k-1)^{m-1} (k-2)!}{m!} \ \tilde S_{\bar i_1,
\dots, \bar i_{\bar k}} (m, k)\ \leqslant$$
$$\leqslant \ \sum\limits_{m \geqslant 1}
\frac{I^{m-1} (k-1)^m (k-2)!}{m!} \
\tilde S_{\bar i_1, \dots, \bar i_{\bar k}} (m, k) \ \leqslant$$
$$\leqslant \ \sum\limits_{m \geqslant 1}\frac{I^{m-1} (k-1)^m (k-2)!}{m!}
m (\bar k-1)!  C_{\bar I \bar k - \bar k}^{k-2}
C_{\bar I \bar k}^{\bar k-m} \ =$$
$$=\ (k-1)! (\bar k-1)! \ \sum\limits_{m \geqslant 1}\frac{I^{m-1}
(k-1)^{m-1}}{(m-1)!}\ C_{\bar I \bar k - \bar k}^{k-2}
\ C_{\bar I \bar k}^{\bar k-m}\ \leqslant$$
$$\leqslant \ (k-1)!  (\bar k-1)!  e^{I(k-1)}  2^{\bar I \bar k - \bar k}
2^{\bar I \bar k}\ . \qed$$
\enddemo

\demo{Proof of Theorem 4.2}
The coefficient of
$t_0^{\dots} t_{i_1}^{n_1} \dots t_{i_I}^{n_I} \bar
t_{\bar i_1}^{\bar n_1} \dots \bar t_{\bar i_{\bar I}}^{\bar n_{\bar I}}$
is equal to
$$\frac{i_1^{n_1}  \dots  i_I^{n_I}}{n_1!  \dots  n_I!} \
\frac{\bar i_1^{\bar n_1}  \dots  \bar i_{\bar I}^{\bar n_{\bar I}}}
{\bar n_1!  \dots  \bar n_{\bar I}!}
\ N_i^2 \left(\matrix i_1,&\dots,&i_I \\ n_1,&\dots,&n_I \endmatrix \left|
\matrix \bar i_1,&\dots,&\bar i_{\bar I} \\ \bar n_1,&\dots,
&\bar n_{\bar I} \endmatrix \right.\right) \ =$$
$$= \ \frac{i_1^{n_1}  \dots  i_I^{n_I}}{n_1!  \dots 
n_I!} \ \frac{\bar i_1^{\bar n_1}  \dots 
\bar i_{\bar I}^{\bar n_{\bar I}}}{\bar n_1!  \dots  \bar n_{\bar I}!}
\ N_i^1 \left(\underbrace{i_1, \dots, i_1}_{n_1}, \dots,
\underbrace{i_I, \dots, i_I}_{n_I} \left|
\underbrace{\bar i_1, \dots, \bar i_1}_{\bar n_1}, \dots,
\underbrace{\bar i_{\bar k}, \dots,
\bar i_{\bar I}}_{\bar n_{\bar I}} \right.\right)\ \leqslant$$

$$\leqslant \ \frac{i_1^{n_1}  \dots  i_I^{n_I}}{n_1!  \dots  n_I!} \
\frac{\bar i_1^{\bar n_1}  \dots 
\bar i_{\bar I}^{\bar n_{\bar I}}}{\bar n_1!  \dots 
\bar n_{\bar I}!} 
k!  \bar k!  e^{\tilde I (k-1)}  2^{2\tilde I \bar k} \ \leqslant
\ \tilde I^K  e^{\tilde I K}  2^{\tilde I K}  \frac{k!  \bar k!}{n_1! 
\dots  n_I!  \bar n_1!  \dots  \bar n_{\bar I}!}\ \leqslant $$
$$\ \leqslant \ \tilde I^K  e^{\tilde I K}  2^{\tilde I K}
\tilde I^K \ \leqslant \ (\tilde I^2 2^{\tilde I} e^{\tilde I})^K\ , $$
where
$k=n_1+\dots +n_I$, $\bar k=\bar n_1+\dots +\bar n_I$, $K=k+\bar k$ and 
$\tilde I= max(I,\bar I)$.

Now consider monomials of degree
$K$ in $t_0,t_1,\bar t_1,\dots ,t_n,\bar t_n$ . The number of such monomial is $(2n)^K$. Thus, their sum in the series is at most
$$(n^2  2^n  e^n)^K  (2n)^K  (4  n^3  2^n e^n)^{-K} \ \leqslant 2^{-K}\ .$$
This implies the convergence of the series $F(\tilde t)$.
\qed
\enddemo

\vskip 0.4cm

\centerline {References}
1. A.N.Varchenko, P.I.Etingof, Why the boundary of a round drop becomes 
a curve of order four. University Lecture Series, 3. American Mathematical Society, Providence, RI,1992.

2. P.S.Novikov, On uniqueness of the solution  of the inverse
problem for a potential. Soviet Math.Dok. 1938, v.18 N 3, 165-168.

3. A.Marshakov, A.Zabrodin, On the Dirichlet boundary
problem and Hirota equations, e-print archiv: hep-th/0305259.

4. P.Wiegmann, A.Zabrodin, Conformal maps and dispersionless
integrable hierarchies, Commun. Math.Phys. 213 (2000), 523.

5. K.Takasaki, T.Takebe, Integrable hierarchies and dispersionless
limit. Rev. Math. Phys. 7(1995), 743-808.

6. R.Dijkgraaf, G.Moor, R.Plessner, The partition function
of two-dimensional string theory, Nucl. Phys, B394 (1993), 356-382.

7. L.L.Chau, A.Zaboronsky, On the structure of correlation functions
in the normal matrix model, Commun. Math. Thys., 196, 203 (1998).

8. S.Natanzon, Integrable systems and effectivisation of Riemann
theorem about domains of the complex plane,
e-print archiv: math.CV/0103136.

9. Yu.Klimov, A.Korzh, S.Natanzon, From $2D$ Toda hierarchy
to conformal maps for domains of Riemann sphere, e-print
archiv: math.NA/0212361.

10. A.Hurwitz, R.Courant, Vorlesungen \"uber allgemeine
Funktionentheorie und ellip\-tische Funktionen. Herausgegeben und 
erg\"anzt durch einen Abschnitt \"uber geometrische Funktionentheorie,
Springer-Verlag, 1964.

11. I.K.Kostov, I.Krichever, M.Mineev--Weinstein, P.B.Wiegmann, A.Zabrodin,
\linebreak $\tau$-function for analytic curves, Random matrices and their 
applications, MSRI publication, v.40, Cambridge Acabemic Press, 2001. 

12. S.M.Natanzon, Formulas for $A_n$ and $B_n$--solutions of
WDVV equations, Journal of Geometry and Physics,  39/4, (2001), 
323-336.

13. R.Graham, D.Knuth, O.Patashnik, Concrete mathematics,
Addison-Wesley pub\-li\-shing Company Reading MA, 1994.

14. L.A.Takhtajan, Free boson tau-functions for compact Riemann surfaces
and closed smooth Jordan curves. Currents correlation functions.
Lett. Math. Phys. 56 (2001), 181-228.

15. A.Marshakov, P.B.Wiegmann, A.Zabrodin, Integrable structure of the
Dirichlet boundary problem in two dimentions, Commun.Math.Phys.227(2002)
131-153

16. I.Krichever, A.Marshakov, A.Zabrodin, Integrable structure of the
Dirichlet boundary problem in multiply-connected domains,
e-print archiv: hep-th/0309010.

\end